\documentclass[a4paper,11pt]{amsart}
\usepackage{graphicx}
\usepackage{mathptmx}
\usepackage{amsmath}
\usepackage{amssymb}
\usepackage{enumitem}
\usepackage{xcolor}

\newmuskip\pFqmuskip

\newcommand*\pFq[6][8]{%
  \begingroup 
  \pFqmuskip=#1mu\relax
  \mathcode`=\string"8000
  \begingroup\lccode`\~=`\,
  \lowercase{\endgroup\let~}\pFqcomma
  F^{#2}_{#3}{\left(\genfrac..{0pt}{}{#4}{#5}\bigg|#6\right)}%
  \endgroup
}
\newcommand{\pFqcomma}{\mskip\pFqmuskip}

\newtheorem{theorem}{Theorem}
\newtheorem{lemma}[theorem]{Lemma}
\newtheorem{corollary}[theorem]{Corollary}

\newtheorem{remark}[theorem]{Remark}

\begin{document}

\title[poly-Dedekind type $DC$ sums involving poly-Euler functions]{poly-Dedekind type $DC$ sums involving poly-Euler functions}

\author{Yuankui Ma}
\address{School of Science, Xi’an Technological University, Xi’an, 710021, Shaanxi, P. R. China}
\email{mayuankui@xatu.edu.cn}

\author{Dae San Kim}
\address{Department of Mathematics, Sogang University, Seoul 121-742, Republic of Korea}
\email{dskim@sogang.ac.kr}

\author{Hyunseok Lee}
\address{Department of Mathematics, Kwangwoon University, Seoul 139-701, Republic of Korea}
\email{luciasconstant@kw.ac.kr}

\author{Hanyoung Kim}
\address{Department of Mathematics, Kwangwoon University, Seoul 139-701, Republic of Korea}
\email{gksaud213@kw.ac.kr}

\author{Taekyun  Kim}
\address{School of Science, Xi’an Technological University, Xi’an, 710021, Shaanxi, P. R. China\\
Department of Mathematics, Kwangwoon University, Seoul 139-701, Republic of Korea}
\email{tkkim@kw.ac.kr}


\maketitle

\begin{abstract}
The classical Dedekind sums appear in the transformation behavior of the logarithm of
the Dedekind eta-function under substitutions from the modular group.
The Dedekind sums and their generalizations are defined in terms of Bernoulli
functions and their generalizations, and are shown to satisfy some reciprocity relations.
In contrast, Dedekind type DC (Daehee and Changhee) sums and their generalizations are defined in terms
of Euler functions and their generalizations. The purpose of this paper is to introduce
the poly-Dedekind type DC sums, which are obtained from the Dedekind type DC
sums by replacing the Euler function by poly-Euler functions of arbitrary indices,
and to show that those sums satisfy, among other things, a reciprocity relation.
\end{abstract}

\section{Introduction}
It is well known that Euler polynomials are defined by
\begin{equation}
\frac{2}{e^{t}+1}e^{xt}=\sum_{n=0}^{\infty}E_{n}(x)\frac{t^{n}}{n!},\quad(\mathrm{see}\
[1-3,5-7,10-17,19-21]). \label{1}
\end{equation}
When $x=0$, $E_{n}=E_{n}(0)$ are called the Euler numbers. \par
From \eqref{1}, we note that
\begin{equation}
E_{n}(x)=\sum_{l=0}^{n}\binom{n}{l}E_{l}x^{n-l},\quad(n\ge
0),\quad(\mathrm{see}\ [1-3,5-7,10-17,19-22]). \label{2}
\end{equation}
The first few of Euler numbers are $E_0=1,E_1=-\frac{1}{2},E_2=0,E_3=\frac{1}{4},E_4=0,E_5=-\frac{1}{2},\dots$,
and $E_{2k}=0$, for $k=1,2,\dots$. \par
From \eqref{1}, we note that $E_{0}=1$, $E_{n}(1)+E_{n}=2\delta_{0,n}$, $(n\ge 0)$, where $\delta_{n,k}$ is
the Kroneceker's symbol. The Euler functions $\overline{E}_{n}(x)$ are defined by
\begin{equation}
\overline{E}_{n}(x)=E_{n}\big(x-[x]\big),\quad(n\ge
0),\quad(\mathrm{see}\ [1,6,14,20]),\label{3}
\end{equation}
where $[x]$ denotes the greatest integer not exceeding $x$. \par
From \eqref{1}, we can easily derive the following identity
\begin{equation}
2\sum_{k=0}^{n-1}(-1)^{k}k^{l}=(-1)^{n-1}E_{l}(n)+E_{l},\quad(n\in\mathbb{N}). \label{4}
\end{equation}
It is known that Dedekind type $DC$ sums are given by
\begin{equation}
T_{p}(h,m)=2\sum_{\mu=0}^{m-1}(-1)^{\mu}\frac{\mu}{m}\overline{E}_{p}\bigg(\frac{h\mu}{m}\bigg),
\quad(h,m\in\mathbb{N}),\quad(\mathrm{see}\ [14,20]).\label{5}
\end{equation}
Note that
\begin{displaymath}
T_{1}(h,m)=
2\sum_{\mu=0}^{m-1}(-1)^{\mu}\bigg(\!\bigg(\frac{\mu}{m}\bigg)\!\bigg)
\bigg(\!\bigg(\frac{h\mu}{m}\bigg)\!\bigg),\quad(\mathrm{see}\
[1,2,6,12,13,21].
\end{displaymath}
where $(\!(x)\!)$ is defined by
\begin{displaymath}
(\!(x)\!)=\left\{\begin{array}{ccc}
x-[x]-\frac{1}{2}, & \textrm{if $x$ is not an integer,}\\
0, & \textrm{if $x$ is an integer.}\quad\
\end{array}\right.
\end{displaymath}
The Genocchi polynomials are defined by
\begin{equation}
\frac{2t}{e^{t}-1}e^{xt}=\sum_{n=0}^{\infty}G_{n}(x)\frac{t^{n}}{n!},\quad(\mathrm{see}\
[7,11,17]).\label{6}
\end{equation}
When $x=0$, $G_{n}=G_{n}(0)$ are called the Genocchi numbers. \par
Note that $G_{0}=0$, $G_{1}=1,\ G_{2}=-1,\ G_{3}=0,\ G_{4}=1,\ G_{5}=0, \ G_{6}=-3,
\dots,$ and $G_{2k+1}=0$ for $k=1,2,3,\dots.$\par
By \eqref{1} and \eqref{6}, we get
\begin{displaymath}
\frac{G_{n+1}(x)}{n+1}=E_{n}(x),\quad \frac{G_{n+1}}{n+1}=E_{n},\quad (n\ge 0).
\end{displaymath}
The degenerate Hardy's polyexponential function of index $k$ is defined by
\begin{equation}
\mathrm{Ei}_{k,\lambda}(x)=\sum_{n=1}^{\infty}\frac{x^{n}(1)_{n,\lambda}}{n^{k}(n-1)!},
\quad(k\in\mathbb{Z}),\quad(\mathrm{see}\ [15]), \label{7}
\end{equation}
where $(x)_{0,\lambda}=1$, $(x)_{n,\lambda}= x(x-\lambda)\cdots\big(x-(n-1)\lambda\big)$, $(n\ge 1)$. \par
Recently, the degenerate poly-Genocchi polynomials of index $k$ are defined in terms of the degenerate Hardy's polyexponential function of index $k$ by
\begin{equation}
\frac{2\mathrm{Ei}_{k,\lambda}\big(\log_{\lambda}(1+t)\big)}{e_{\lambda}(t)+1}e_{\lambda}^{x}(t)
=\sum_{n=0}^{\infty}G_{n,\lambda}^{(k)}(x)\frac{t^{n}}{n!},\quad(\mathrm{see}\
[17]), \label{8}
\end{equation}
where
\begin{displaymath}
e_{\lambda}^{x}(t)=\sum_{n=0}^{\infty}\frac{(x)_{n,\lambda}}{n!}t^{n},\quad e_{\lambda}(t)=e_{\lambda}^{1}(t),
\quad\mathrm{and}\quad \log_{\lambda}(t)=\frac{1}{\lambda}(t^{\lambda}-1)
\end{displaymath}
is the compositional inverse to $e_{\lambda}(t)$ satisfying $e_{\lambda}(\log_{\lambda}(t))
=\log_{\lambda}(e_{\lambda}(t))=t$.
Taking $\lambda\rightarrow 0$ in \eqref{8}, we have the poly-Genocchi polynomials of index $k$ given by
\begin{equation}
\frac{2\mathrm{Ei}_{k}\big(\log(1+t)\big)}{e^{t}+1}e^{xt}
=\sum_{n=0}^{\infty}G_{n}^{(k)}(x)\frac{t^{n}}{n!},\quad(\mathrm{see}\
[7,17]),\label{9}
\end{equation}
where $G_{n}^{(k)}(x)=\lim_{\lambda\rightarrow 0}G_{n,\lambda}^{(k)}(x),$ $(n\ge 0)$, and
\begin{equation}
\mathrm{Ei}_{k}(x)=\sum_{n=1}^{\infty}\frac{x^{n}}{n^{k}(n-1)!},\quad
(\mathrm{see}\ [10,15]).\label{9-1}
\end{equation}
is the polyexponential function of index $k$. \par
\indent When $x=0$, $G_{n}^{(k)}=G_{n}^{(k)}(0)$, $(n\ge 0)$, are called the poly-Genocchi numbers of index $k$.
By \eqref{9}, we easily get $G_{0}^{(k)}=0,\ G_{1}^{(k)}=1,\ G_{2}^{(k)}=-2+2^{1-k},\dots$.
Also, from \eqref{9}, we note that
\begin{equation}
G_{n}^{(k)}(x)=\sum_{l=0}^{n}\binom{n}{l}G_{l}^{(k)}x^{n-l},\quad(n\ge 0). \label{10}
\end{equation} \par

\begin{remark}
The polyexponential functions were first considered by Hardy, which are given by
\begin{equation*}
e(x,\
a|s)=\sum_{n=0}^{\infty}\frac{x^{n}}{(n+a)^{s}n!},\quad(\mathrm{Re}(a)>0),\quad
(\mathrm{see}\ [4,8,9]).
\end{equation*}
In [18], Komatsu defined the polylogarithm factorial function
$Lif_{k}(x)$ by $xLif_{k}(x)=xe(x,1|k)=Ei_{k}(x)$. So the
polylogarithm factorial functions are special cases of Hardy's
polyexponential functions, but our polyexponential functions are
not. In fact, the slight difference between ours and Komatsu's is
crutial in defining, for example, the type 2 poly-Bernoulli
polynomials (see [10,15]) and also in constructing poly-Dedekind
sums associated with such polynomials (see [16,19]). Here we
recall from [10] that the type 2 poly-Bernoulli polynomials
$\beta_n^{(k)}(x)$ of index $k$ are defined by
\begin{equation}
\frac{Ei_k(\log(1+t))}{e^t -1}e^{xt}=\sum_{n=0}^{\infty}\beta_n^{(k)}(x)\frac{t^n}{n!}.\label{9-1}
\end{equation}
We also recall that, for any integer $k$, the poly-Bernoulli polynomials $B_n^{(k)}(x)$ of index $k$ are defined by
\begin{equation}
\frac{Li_k(1-e^{-t})}{1-e^{-t}}e^{xt}=\sum_{n=0}^{\infty}B_n^{(k)}(x)\frac{t^n}{n!}, \label{9-2}
\end{equation}
where the polylogarithm functions $Li_k(x)$ are given by
$Li_k(x)=\sum_{n=1}^{\infty}\frac{x^n}{n^k}$.\par The reason why
$Ei_k(x)$ is needed and $Lif_k(x)$ is not  in \eqref{9-1} is
twofold. The first reason is that $Ei_k(x)$ has order 1, so that
the composition $Ei_k(\log(1+t))$ still has order 1, which is
definitely required, whereas $Lif_{k}(x)$ has order 0, so that
$Lif_k(\log(1+t))$ also has order 0. The second reason is that we
want $\beta_n^{(1)}(x)$ to be the ordinary Bernoulli polynomials
when $k=1$. Indeed, $Ei_1(x)=e^x-1$, so that  $\beta_n^{(1)}(x)$
are those polynomials with $\beta_1^{(1)}(x)=x-\frac{1}{2}$. \par
The construction of the type 2 poly-Bernoulli polynomials are in
parallel with that of the poly-Bernoulli polynomials. We note
$Li_k(x)$ has order 1, so that the composition $Li_k(1-e^{-t})$
also has order 1. In addition, $Li_1(x)=-\log(1-x)$, so that
$B_n^{(1)}(x)$ are the ordinary Bernoulli polynomials with
$B_1^{(1)}(x)=x+\frac{1}{2}$ (see \eqref{9-2}). Thus we may say
that $Ei_k(x)$ is a kind of a compositional inverse type to
$Li_k(x)$.
\end{remark}

Now, we define the {\it{poly-Euler polynomials of index $k$}} by
\begin{equation}
E_{n}^{(k)}(x)=\frac{G_{n+1}^{(k)}(x)}{n+1},\quad(n\ge 0).\label{11}
\end{equation}
When $x=0$, $E_{n}^{(k)}=E_{n}^{(k)}(0)$ are called the poly-Euler numbers of index $k$.
Note that $E_{n}^{(1)}(x)=E_{n}(x)$ and $G_{n}^{(1)}(x)=G_{n}(x)$. \par
From \eqref{10}, we note that
\begin{align}
    E_{n}^{(k)}(x)\ &=\ \frac{1}{n+1} G_{n+1}^{(k)}(x)\ =\ \frac{1}{n+1}\sum_{l=0}^{n+1}\binom{n+1}{l}G_{l}^{(k)}x^{n+1-l} \label{12} \\
    &=\ \frac{1}{n+1}\sum_{l=1}^{n+1}\binom{n+1}{l}G_{l}^{(k)}x^{n+1-l}\ =\ \frac{1}{n+1}\sum_{l=0}^{n}\binom{n+1}{l+1}G_{l+1}^{(k)}x^{n-l} \nonumber\\
    &=\ \sum_{l=0}^{n}\binom{n}{l}\frac{G_{l+1}^{(k)}}{l+1}x^{n-l}\ =\ \sum_{l=0}^{n}\binom{n}{l}E_{l}^{(k)}x^{n-l},\quad(n\ge 0). \nonumber
\end{align}
Apostol considered the generalized Dedekind sums which are given by
\begin{equation}\label{12-1}
S_p(h,m)= \sum_{\mu=1}^{m-1}\frac{\mu}{m}\overline{B}_{p}\big(\frac{h\mu}{m}\big),
\end{equation}
and showed that they satisfy a reciprocity relation in [1,2].
Here $\overline{B}_p(x)=B_p(x-[x])$ are the Bernoulli functions with $B_p(x)$ the Bernoulli polynomials given by
\begin{equation*}
\frac{t}{e^t-1}e^{xt}=\sum_{p=0}^{\infty}B_p(x)\frac{t^p}{p!}.
\end{equation*}
We remark that the Dedekind sum $S(h,m)=S_1(h,m)$ appears in the transformation behaviour of the logarithm of
the Dedekind eta--function under substitutions from the modular group and a reciprocity law of that
was demonstrated by Dedekind in 1892. \par
As an extension of the sums in \eqref{12-1}, the poly-Dedekind sums, which are given by
\begin{equation*}
S_p^{(k)}(h,m)= \sum_{\mu=1}^{m-1}\frac{\mu}{m}\overline{B}_{p}^{(k)}\bigg(\frac{h\mu}{m}\bigg),
\end{equation*}
were considered and a reciprocity law for those sums was shown in
[16,19]. Here $B_{p}^{(k)}(x)$ are the type 2 poly-Bernoulli
polynomials of index $k$,
$\overline{B}_{p}^{(k)}(x)=B_{p}^{(k)}(x-[x])$, and
$B_p^{(1)}(x)=B_p(x)$, (see [16]). \par

The Dedekind type DC sums (see \eqref{5}) were first introduced
and shown to satisfy a reciprocity relation in [14]. The aim of
this paper is to introduce the poly-Dedekind type DC sums (see
\eqref{12-2}), which are obtained from the Dedekind type DC sums
by replacing the Euler function by poly-Euler functions of
arbitrary indices, and to show that those sums satisfy, among
other things, a reciprocity relation (see \eqref{12-3}). The
motivation of this paper is to explore our new sums in connection
with modular forms, zeta fuctions and trigonometric sums, just as
in the cases of Apostol-Dedekind sums, their generalizations and
of some related sums. Indeed, Simsek [20] found trigonometric
representations of the Dedekind type DC sums and their relations
to  Clausen functions, polylogarithm function, Hurwitz zeta
function, generalized Lambert series (G-series), and Hardy-Berndt
sums. In addition, Bayad and Simsek [3] studied three new shifted
sums of Apostol-Dedekind-Rademacher type. These sums generalize
the classical Dedekind-Rademacher sums and can be expressed in
terms of Jacobi modular forms or cotangent functions or special
values of the Barnes multiple zeta functions. They found
reciprocity laws for these sums and demonstrated that some
well-known reciprocity laws can be deduced from their results. We
plan to carry out this line of research in a subsequent paper.\par
In this paper, we consider the poly-Dedekind type DC sums defined
by
\begin{equation}
T_{p}^{(k)}(h,m)=2\sum_{\mu=1}^{m-1}(-1)^{\mu}\frac{\mu}{m}\overline{E}_{p}^{(k)}\bigg(\frac{h\mu}{m}\bigg),\label{12-2}
\end{equation}
where $h,m,p \in \mathbb{N}$, and $\overline{E}_{p}^{(k)}$ are the poly-Euler functions of index $k$ given by $\overline{E}_{p}^{(k)}(x)=E_{p}^{(k)}\big(x-[x]\big)$.
We show the following reciprocity relation for the poly-Dedekind type DC sums given by
\begin{align}
&m^{p}T_{p}^{(k)}(h,m)+h^{p}T_{p}^{(k)}(m,h)\label{12-3}\\
& =2\sum_{\mu=0}^{m-1}\sum_{l=0}^{p}\sum_{\nu=0}^{h-1}\sum_{j=1}^{p+1-l}(-1)^{\mu+\nu}\frac{(mh)^{l-1}\binom{p}{l}S_{1}(p-l+1,j)}{(p-l+1)j^{k-1}}\big((\mu h)m^{p-l}+(\nu m)h^{p-l}\big)\overline{E}_{l}\bigg(\frac{\nu}{h}+\frac{\mu}{m}\bigg)\nonumber,
\end{align}
where $m,h,p\in\mathbb{N}$ with $m\equiv 1$ $(\mathrm{mod}\ 2)$, $h\equiv 1$ $(\mathrm{mod}\ 2)$, and $k\in\mathbb{Z}$. \par
For $k=1$, this reciprocity relation for the poly-Dedekind type DC sums reduces to that for the Dedekind type DC sums given by (see Corollary 15)
\begin{align*}
&m^{p}T_{p}(h,m)+h^{p}T_{p}(m,h) \\
&=2(mh)^{p-1}\sum_{\mu=0}^{m-1}\sum_{\nu=0}^{h-1}(-1)^{\mu+\nu-1}(\mu h+\nu m)\overline{E}_{p}\bigg(\frac{\nu}{h}+\frac{\mu}{m}\bigg),
\end{align*}
where $m,h,p\in\mathbb{N}$ with $m\equiv 1$ $(\mathrm{mod}\ 2)$, $h\equiv 1$ $(\mathrm{mod}\ 2)$. \par
\begin{remark}
From \eqref{5} and \eqref{12-2}, we see that the poly-Dedekind type DC sums are obtained from the Dedekind type DC sums by replacing the Euler functions by poly-Euler functions of arbitrary indices. Here we have to observe that the key to this generalization is the construction of poly-Euler polynomials defined in \eqref{11}, which is done in an elaborate manner. First, we replace $t$ by $Ei_k\big(\log(1+t)\big)$ as in \eqref{9}, so that we construct the poly-Genocchi polynomials $G_n^{(k)}(x)$ of index $k$ such that $G_n^{(1)}(x)=G_n(x)$ are the usual Genocchi polynomials. Next, by defining the poly-Euler polynomials $E_n^{(k)}(x)$ as in \eqref{11}, so that we have the desirable property $E_n^{(1)}(x)=E_n(x)$. Consequently, for $k=1$, the poly-Dedekind type DC sums $T_p^{(k)}(h,m)$ in \eqref{12-2} reduces to the Dedekind type DC sums $T_p(h,m)$ in \eqref{5}.
\end{remark}
In Section 2, we will derive various facts about the poly-Genocchi polynomials and poly-Euler polynomials that will be needed in the next section. In Section 3, we will define the poly-Dedekind type DC sums and demonstrate, among other things, a reciprocity relation for them.

\section{Poly-Genocchi polynomials and Poly-Euler polynomials}
By \eqref{9}, we have
 \begin{align}
 2\mathrm{Ei}_{k}\big(\log(1+t)\big)\ &=\ \frac{2\mathrm{Ei}_{k}\big(\log(1+t)\big)}{e^t+1}e^t+ \frac{ 2\mathrm{Ei}_{k}\big(\log(1+t)\big)}{e^t+1}\label{13} \\
 &=\ \sum_{n=0}^{\infty}\big(G_{n}^{(k)}(1)+G_{n}^{(k)}\big)\frac{t^{n}}{n!}. \nonumber
\end{align}
On the other hand, we also have
\begin{align}
 2\mathrm{Ei}_{k}\big(\log(1+t)\big)\ &=\ 2\sum_{m=1}^{\infty}\frac{1}{m^{k}(m-1)!}\big(\log(1+t)\big)^{m} \label{14} \\
 &=\ \sum_{n=1}^{\infty}\bigg(2\sum_{m=1}^{n}\frac{1}{m^{k-1}}S_{1}(n,m)\bigg)\frac{t^{n}}{n!},\nonumber
\end{align}
where $S_{1}(n,m)$ are the Stirling numbers of the first kind. \par
Therefore, by \eqref{13} and \eqref{14}, we get the following theorem.
\begin{theorem}
For $n\ge 1$, we have
\begin{displaymath}
2\sum_{m=1}^{n}\frac{1}{m^{k-1}}S_{1}(n,m)=G_{n}^{(k)}(1)+G_{n}^{(k)}.
\end{displaymath}
\end{theorem}
\begin{corollary}
For $n\ge 1$, we have
\begin{displaymath}
\frac{2}{n}\sum_{m=1}^{n}\frac{1}{m^{k-1}}S_{1}(n,m)=E_{n-1}^{(k)}(1)+E_{n-1}^{(k)}.
\end{displaymath}
\end{corollary}
From \eqref{10} and \eqref{12}, we see that
\begin{displaymath}
\frac{d}{dx}G_{n+1}^{(k)}(x)=(n+1)G_{n}^{(k)}(x),\quad\frac{d}{dx}E_{n}^{(k)}(x)=nE_{n-1}^{(k)}(x),\quad(n\ge 1).
\end{displaymath}
Thus, we note that
\begin{align*}
\int_{0}^{x}G_{n}^{(k)}(x)dx\ &=\ \frac{1}{n+1}\big(G_{n+1}^{(k)}(x)-G_{n+1}^{(k)}\big), \\
\int_{0}^{x}E_{n-1}^{(k)}(x)dx\ &=\ \frac{1}{n}\big(E_{n}^{(k)}(x)-E_{n}^{(k)}\big),\quad(n\ge 1).
\end{align*}
From \eqref{2} and \eqref{8}, we have
\begin{align}
\frac{2\mathrm{Ei}_{k}\big(\log(1+t)\big)}{e^{t}+1}e^{xt}\ &=\ \sum_{m=1}^{\infty}\frac{\big(\log(1+t)\big)^{m}}{m^{k}(m-1)!}\sum_{l=0}^{\infty}E_{l}(x)\frac{t^{l}}{l!} \label{15} \\
&=\ \sum_{j=1}^{\infty}\sum_{m=1}^{j}\frac{S_{1}(j,m)}{m^{k-1}}\frac{t^{j}}{j!}\sum_{l=0}^{\infty}E_{l}(x)\frac{t^{l}}{l!} \nonumber\\
&=\ \sum_{n=1}^{\infty}\bigg(\sum_{j=1}^{n}\sum_{m=1}^{j}\binom{n}{j}\frac{S_{1}(j,m)}{m^{k-1}}E_{n-j}(x)\bigg)\frac{t^{n}}{n!}. \nonumber
\end{align}
On the other hand, we also have
\begin{equation}
\frac{2\mathrm{Ei}_{k}\big(\log(1+t)\big)}{e^{t}+1}e^{xt}=\sum_{n=1}^{\infty}G_{n}^{(k)}(x)\frac{t^{n}}{n!}=\sum_{n=1}^{\infty}nE_{n-1}^{(k)}(x)\frac{t^{n}}{n!}. \label{16}
\end{equation}
Therefore, by \eqref{15} and \eqref{16}, we obtain the following theorem.
\begin{theorem}
For $n\in\mathbb{N}$, we have
\begin{displaymath}
E_{n-1}^{(k)}(x)=\frac{1}{n}\sum_{j=1}^{n}\sum_{m=1}^{j}\binom{n}{j}\frac{S_{1}(j,m)}{m^{k-1}}E_{n-j}(x).
\end{displaymath}
\end{theorem}
For $m\in\mathbb{N}$ with $m\equiv 1$ $(\mathrm{mod}\ 2)$, we have
\begin{align}
\frac{2}{e^{t}+1}e^{xt}\ &=\ \frac{2}{1+e^{mt}}\sum_{i=0}^{m-1}(-1)^{i}e^{it}e^{xt} \label{17}\\
&=\ \sum_{i=0}^{m-1}(-1)^{i}\frac{2}{e^{mt}+1}e^{\big(\frac{i+x}{m}\big)mt}\nonumber \\
&=\ \sum_{n=0}^{\infty}\bigg(m^{n}\sum_{i=0}^{m-1}(-1)^{i}E_{n}\bigg(\frac{x+i}{m}\bigg)\bigg)\frac{t^{n}}{n!}.\nonumber
\end{align}
By \eqref{1} and \eqref{17}, we get the distribution relation
\begin{equation}
E_{n}(x)=m^{n}\sum_{i=0}^{m-1}(-1)^{i}E_{n}\bigg(\frac{x+i}{m}\bigg), \label{18}
\end{equation}
where $m\in\mathbb{N}$ with $m\equiv 1$ $(\mathrm{mod}\ 2)$, and $n\ge 0$. \par
For $x\in\mathbb{N}$, we have
\begin{align}
&2\sum_{i=0}^{x-1}(-1)^{i}e^{it}\mathrm{Ei}_{k}\big(\log(1+t)\big)\ =\ 2\sum_{i=0}^{x-1}(-1)^{i}\sum_{l=0}^{\infty}i^{l}\frac{t^{l}}{l!} \sum_{j=1}^{\infty}\frac{\big(\log(1+t)\big)^{j}}{j^{k}(j-1)!}\label{19}\\
&=\ \sum_{l=0}^{\infty}2\sum_{i=0}^{x-1}(-1)^{i}i^{l}\frac{t^{l}}{l!}\sum_{m=1}^{\infty}
\sum_{j=1}^{m}\frac{S_{1}(m,j)}{j^{k-1}}\frac{t^{m}}{m!}\nonumber \\
&=\ \sum_{n=1}^{\infty}\bigg(2\sum_{m=1}^{n}\sum_{j=1}^{m}\sum_{i=0}^{x-1}(-1)^{i}i^{n-m}\binom{n}{m}\frac{S_{1}(m,j)}{j^{k-1}}\bigg)\frac{t^{n}}{n!}.\nonumber
\end{align}
On the other hand, we also have
\begin{align}
2\sum_{i=0}^{x-1}(-1)^{i}e^{it}\mathrm{Ei}_{k}\big(\log(1+t)\big)\ &=\ \frac{2\mathrm{Ei}_{k}\big(\log(1+t)\big)}{e^{t}+1}\big((-1)^{x-1}e^{xt}+1\big) \label{20} \\
&=\ \sum_{n=0}^{\infty}\big((-1)^{x-1}G_{n}^{(k)}(x)+G_{n}^{(k)}\big)\frac{t^{n}}{n!}.\nonumber
\end{align}
Therefore, by \eqref{19} and \eqref{20}, we obtain the following theorem.
\begin{theorem}
For $x, n \in\mathbb{N}$, we have
\begin{displaymath}
(-1)^{x-1}G_{n}^{(k)}(x)+G_{n}^{(k)}=2\sum_{m=1}^{n}\sum_{j=1}^{m}\sum_{i=0}^{x-1}(-1)^{i}i^{n-m}\binom{n}{m}\frac{S_{1}(m,j)}{j^{k-1}}.
\end{displaymath}
\end{theorem}
Note that, for $k=1$, we have
\begin{displaymath}
(-1)^{x-1}G_{n}^{(1)}(x)+G_{n}^{(1)}=2n\sum_{i=0}^{x-1}(-1)^{i}i^{n-1}.
\end{displaymath}
\begin{corollary}
For $x,n\in\mathbb{N}$, we have
\begin{displaymath}
(-1)^{x-1}E_{n-1}^{(k)}(x)+E_{n-1}^{(k)}=\frac{2}{n}\sum_{m=1}^{n}\sum_{j=1}^{m}\sum_{i=0}^{x-1}(-1)^{i}i^{n-m}\binom{n}{m}\frac{S_{1}(m,j)}{j^{k-1}}.
\end{displaymath}
\end{corollary}
Note that
\begin{displaymath}
(-1)^{x-1}E_{n-1}^{(1)}(x)+E_{n-1}^{(1)}=2\sum_{i=0}^{x-1}(-1)^{i}i^{n-1},\quad (n,x\in\mathbb{N}).
\end{displaymath}

For $m\in\mathbb{N}$ with $m\equiv 1\ (\mathrm{mod}\ 2)$, we note that
\begin{align}
&2\frac{\mathrm{Ei}_{k}\big(\log(1+t)\big)}{e^{t}+1}e^{xt}\ =\ \frac{1}{m}\sum_{s=0}^{m-1}(-1)^{s}e^{\big(\frac{s+x}{m}\big)mt}\frac{2mt}{e^{mt}+1}\frac{1}{t}\mathrm{Ei}_{k}\big(\log(1+t)\big)\label{21} \\
&=\ \sum_{l=0}^{\infty}m^{l-1}\sum_{s=0}^{m-1}(-1)^{s}G_{l}\bigg(\frac{s+x}{m}\bigg)\frac{t^{l}}{l!} \frac{1}{t}\sum_{j=1}^{\infty}\frac{\big(\log(1+t)\big)^{j}}{j^{k}(j-1)!}\nonumber
\end{align}
\begin{align*}
&=\ \sum_{l=0}^{\infty}m^{l-1}\sum_{s=0}^{m-1}(-1)^{s}G_{l}\bigg(\frac{s+x}{m}\bigg)\frac{t^{l}}{l!}\frac{1}{t}\sum_{i=1}^{\infty}\sum_{j=1}^{i}\frac{S_{1}(i,j)}{j^{k-1}}\frac{t^{i}}{i!} \nonumber  \\
&=\ \sum_{l=0}^{\infty}m^{l-1}\sum_{s=0}^{m-1}(-1)^{s}G_{l}\bigg(\frac{s+x}{m}\bigg)\frac{t^{l}}{l!}\sum_{i=0}^{\infty}\sum_{j=1}^{i+1}\frac{1}{j^{k-1}}\frac{S_{1}(i+1,j)}{i+1}\frac{t^{i}}{i!}\nonumber \\
&=\ \sum_{n=0}^{\infty}\bigg(\sum_{l=0}^{n}\binom{n}{l}m^{l-1}\sum_{j=1}^{n-l+1}\sum_{s=0}^{m-1}(-1)^{s}G_{l}\bigg(\frac{s+x}{m}\bigg)\frac{1}{j^{k-1}}\frac{S_{1}(n-l+1,j)}{n-l+1}\bigg)\frac{t^{n}}{n!}.\nonumber
\end{align*}
Therefore, by \eqref{21}, we obtain the following theorem.
\begin{theorem}
For $n \ge 0$, and $m\in\mathbb{N}$ with $m \equiv 1$ $(\mathrm{mod}\ 2)$, we have
\begin{displaymath}
G_{n}^{(k)}(x)=\sum_{l=0}^{n}\binom{n}{l}m^{l-1}\sum_{j=1}^{n-l+1}\sum_{s=0}^{m-1}(-1)^{s}G_{l}\bigg(\frac{s+x}{m}\bigg)\frac{1}{j^{k-1}}\frac{S_{1}(n-l+1,j)}{n-l+1}.
\end{displaymath}
\end{theorem}
From Theorem 6, we have
\begin{align*}
    \frac{G_{n}^{(k)}(x)}{n}\ &=\ \frac{1}{n}\sum_{l=1}^{n}\binom{n}{l}m^{l-1}\sum_{j=1}^{n-l+1}\sum_{s=0}^{m-1}(-1)^{s}G_{l}\bigg(\frac{s+x}{m}\bigg)\frac{1}{j^{k-1}}\frac{S_{1}(n-l+1,j)}{n-l+1}\\
    &=\ \frac{1}{n}\sum_{l=0}^{n-1}\binom{n}{l+1}m^{l}\sum_{j=1}^{n-l}\sum_{s=0}^{m-1}(-1)^{s}G_{l+1}\bigg(\frac{s+x}{m}\bigg)\frac{1}{j^{k-1}}\frac{S_{1}(n-l,j)}{n-l}\\
    &=\ \sum_{l=0}^{n-1}\binom{n-1}{l}m^{l}\sum_{j=1}^{n-l}\sum_{s=0}^{m-1}(-1)^{s}\frac{G_{l+1}\big(\frac{s+x}{m}\big)}{l+1}\frac{1}{j^{k-1}}\frac{S_{1}(n-l,j)}{n-l}\\
    &=\ \sum_{l=0}^{n-1}\binom{n-1}{l}m^{l}\sum_{j=1}^{n-l}\sum_{s=0}^{m-1}(-1)^{s}E_{l}\bigg(\frac{s+x}{m}\bigg)\frac{1}{j^{k-1}}\frac{S_{1}(n-l,j)}{n-l},
\end{align*}
where $n,m\in\mathbb{N}$ with $m\equiv 1$ $(\mathrm{mod}\ 2)$.
Thus we obtain the important corollary that will be used in deriving the reciprocity law in Theorem 14.
\begin{corollary}
For $n, m\in\mathbb{N}$ with $m\equiv 1$ $(\mathrm{mod}\ 2)$, we have
\begin{displaymath}
E_{n-1}^{(k)}(x)=\sum_{l=0}^{n-1}\binom{n-1}{l}m^{l}\sum_{j=1}^{n-l}\sum_{s=0}^{m-1}(-1)^sE_{l}\bigg(\frac{s+x}{m}\bigg)\frac{1}{j^{k-1}}\frac{S_{1}(n-l,j)}{n-l}.
\end{displaymath}
\end{corollary}
Note that
\begin{displaymath}
    E_{n-1}^{(1)}(x)=m^{n-1}\sum_{s=0}^{m-1}(-1)^{s}E_{n-1}\bigg(\frac{s+x}{m}\bigg),
\end{displaymath}
where $n,m\in\mathbb{N}$ with $m\equiv 1$ $(\mathrm{mod}\ 2)$. \par

For $p, s\in\mathbb{N}$ with $s<p$, we have
\begin{equation}
\frac{d^{s}}{dx^{s}}\big(xE_{p}^{(k)}(x)\big)\bigg|_{x=1}=s!\binom{p}{s}E_{p-s}^{(k)}(1)+s!\binom{p}{s-1}E_{p-s+1}^{(k)}(1). \label{22}
\end{equation}
On the other hand, by \eqref{12}, we get
\begin{equation}
\frac{d^{s}}{dx^{s}}\big(xE_{p}^{(k)}(x)\big)\bigg|_{x=1}=s!\sum_{\nu=0}^{p}\binom{p-\nu+1}{s}\binom{p}{\nu}E_{\nu}^{(k)}. \label{23}
\end{equation}
Therefore, by \eqref{22} and \eqref{23}, we obtain the following lemma.
\begin{lemma}
For $p, s\in\mathbb{N}$ with $s<p$, we have
\begin{displaymath}
\sum_{\nu=0}^{p}\binom{p-\nu+1}{s}\binom{p}{\nu}E_{\nu}^{(k)}=\binom{p}{s}E_{p-s}^{(k)}(1)+\binom{p}{s-1}E_{p-s+1}^{(k)}(1).
\end{displaymath}
\end{lemma}
In particular, for $k=1$, and $p, s\in\mathbb{N}$ with $p \equiv 1$ $(\mathrm{mod}\ 2)$, $s\equiv 0$ $(\mathrm{mod}\ 2)$, we have
\begin{displaymath}
\sum_{\nu=0}^{p}\binom{p-\nu+1}{s}\binom{p}{\nu}E_{\nu}=\binom{p}{s}E_{p-s}(1)=-\binom{p}{s}E_{p-s}.
\end{displaymath}
We note that
\begin{align}
\int_{0}^{1}xE_{p}^{(k)}(x)dx\ &=\ \frac{1}{p+1}E_{p+1}^{(k)}(1)-\frac{1}{p+1}\int_{0}^{1}E_{p+1}^{(k)}(x)dx\label{24} \\
&=\ \frac{1}{p+1}E_{p+1}^{(k)}(1)-\frac{1}{(p+1)(p+2)}\big(E_{p+2}^{(k)}(1)-E_{p+2}^{(k)}\big).\nonumber
\end{align}
On the other hand, from \eqref{12}, we have
\begin{align}
\int_{0}^{1}xE_{p}^{(k)}(x)dx\ &=\ \sum_{\nu=0}^{p}\binom{p}{\nu}E_{\nu}^{(k)}\int_{0}^{1}x^{p-\nu+1}dx\label{25} \\
&=\ \sum_{\nu=0}^{p}\binom{p}{\nu}\frac{E_{\nu}^{(k)}}{p-\nu+2}. \nonumber
\end{align}
Therefore, by \eqref{24} and \eqref{25}, we obtain the following lemma.
\begin{lemma}
For $p\in\mathbb{N}$, we have
\begin{displaymath}
\sum_{\nu=0}^{p}\binom{p}{\nu}\frac{E_{\nu}^{(k)}}{p-\nu+2}=\frac{1}{p+1}E_{p+1}^{(k)}(1)-\frac{1}{(p+1)(p+2)}E_{p+2}^{(k)}(1)+\frac{1}{(p+1)(p+2)}E_{p+2}^{(k)}.
\end{displaymath}
\end{lemma}
In particular, for $p\in\mathbb{N}$ with $p\equiv 1$ $(\mathrm{mod}\ 2)$, and $k=1$, we get
\begin{displaymath}
    \sum_{\nu=0}^{p}\binom{p}{\nu}\frac{E_{\nu}}{p-\nu+2}=\frac{2}{(p+1)(p+2)}E_{p+2}.
\end{displaymath}

\section{Poly-Dedekind type $DC$ sums}
The Dedekind type $DC$ sums are defined by
\begin{equation}
T_{p}(h,m)=2\sum_{\mu=1}^{m-1}(-1)^{\mu}\frac{\mu}{m}\overline{E}_{p}\bigg(\frac{h\mu}{m}\bigg),\quad(h,m\in\mathbb{N}),\label{26}
\end{equation}
where $\overline{E}_{p}(x)$ is the $p$-th Euler function (see
[14,20]). \par For $p\in\mathbb{N}$ with $p\equiv 1$
$(\mathrm{mod}\ 2)$, and relative prime positive integers  $m,h$
with $m\equiv 1$ $(\mathrm{mod}\ 2)$, $h\equiv 1$ $(\mathrm{mod}\
2)$, the reciprocity law of $T_{p}(h,m)$ is given by
\begin{align*}
    &m^{p}T_{p}(h,m)+h^{p}T_{p}(m,h)\\
    &=2\sum_{\mu}\bigg(mh \bigg(E+\frac{\mu}{m} \bigg)+m\bigg(E+h-\bigg[\frac{h\mu}{m}\bigg]\bigg)\bigg)^{p}\\
    &\qquad+\big(hE+mE\big)^{p}+(p+2)E_{p},
\end{align*}
where $\mu$ runs over all integers satisfying $0 \le \mu \le m-1$ and $\mu-\bigg[\frac{h\mu}{m}\bigg]\equiv 1\ (\mathrm{mod}\ 2)$, and
\begin{displaymath}
\big(hE+mE\big)^{p}=\sum_{l=0}^{p}\binom{p}{l}h^{l}E_{l}h^{l}m^{p-l}E_{p-l}.
\end{displaymath} \par
For the rest of our discussion, {\it{we assume that $k$ is any integer}}.
In light of \eqref{26}, we define {\it{poly-Dedekind type $DC$ sums}} given by
\begin{equation}
T_{p}^{(k)}(h,m)=2\sum_{\mu=1}^{m-1}\bigg(\frac{\mu}{m}\bigg)(-1)^{\mu}\overline{E}_{p}^{(k)}\bigg(\frac{h\mu}{m}\bigg),\label{27}
\end{equation}
where $h,m,p\in\mathbb{N}$, and $\overline{E}_{p}^{(k)}$ are the poly-Euler functions given by
\begin{displaymath}
    \overline{E}_{p}^{(k)}(x)=E_{p}^{(k)}\big(x-[x]\big).
\end{displaymath} \par
By \eqref{26} and \eqref{27}, we get
\begin{equation}
T_{p}^{(1)}(h,m)=2\sum_{\mu=1}^{m-1}(-1)^{\mu}\bigg(\frac{\mu}{m}\bigg)\overline{E}_{p}\bigg(\frac{\mu}{m}\bigg)=T_{p}(h,m).\label{28}
\end{equation}
Let us take $h=1$. Then we have
\begin{align}
T_{p}^{(k)}(1,m)\ &=\ 2\sum_{\mu=0}^{m-1}(-1)^{\mu}\bigg(\frac{\mu}{m}\bigg)\sum_{\nu=0}^{p}\binom{p}{\nu}\bigg(\frac{\mu}{m}\bigg)^{p-\nu}E_{\nu}^{(k)}    \label{29} \\
&=\ \sum_{\nu=0}^{p}\binom{p}{\nu}E_{\nu}^{(k)}m^{-(p-\nu+1)}\bigg(2\sum_{\mu=1}^{m-1}(-1)^{\mu}\mu^{p-\nu+1}\bigg).\nonumber
\end{align}
Assume that  $m\in\mathbb{N}$ with $m\equiv 1$ $(\mathrm{mod}\ 2)$. Then, by \eqref{4} and \eqref{29}, we get
\begin{align}
T_{p}^{(k)}(1,m)\ &=\ \sum_{\nu=0}^{p}\binom{p}{\nu}E_{\nu}^{(k)}m^{-(p+1-\nu)}\big(E_{p+1-\nu}(m)+E_{p+1-\nu}\big) \label{30} \\
&=\ \sum_{\nu=0}^{p}\binom{p}{\nu}E_{\nu}^{(k)}m^{-(p+1-\nu)}\bigg(\sum_{i=0}^{p+1-\nu}\binom{p+1-\nu}{i}m^{p+1-\nu-i}E_{i}+E_{p+1-\nu}\bigg)\nonumber \\
&=\ \sum_{\nu=0}^{p}\binom{p}{\nu}E_{\nu}^{(k)}m^{-(p+1-\nu)}\sum_{i=0}^{p-\nu}\binom{p-\nu+1}{i}m^{p+1-\nu-i}E_{i}\nonumber \\
&\quad +2\sum_{\nu=0}^{p}\binom{p}{\nu}E_{\nu}^{(k)}E_{p+1-\nu}m^{-(p+1-\nu)}.\nonumber
\end{align}
From \eqref{30}, we have
\begin{equation}
m^{p}T_{p}^{(k)}(1,m)=\sum_{\nu=0}^{p}\binom{p}{\nu}E_{\nu}^{(k)}\sum_{i=0}^{p-\nu}\binom{p-\nu+1}{i}m^{p-i}E_{i}+2\sum_{\nu=0}^{p}\binom{p}{\nu}E_{\nu}^{(k)}E_{p+1-\nu}m^{\nu-1}. \label{31}
\end{equation}
Let us define $S_{p}^{(k)}(1,m)$ as
\begin{equation}
m^{p}T_{p}^{(k)}(1,m)-2\sum_{\nu=0}^{p}\binom{p}{\nu}E_{\nu}^{(k)}E_{p+1-\nu}m^{\nu-1}=S_{p}^{(k)}(1,m),
\label{32}
\end{equation}
where $p,m\in\mathbb{N}$ with $m\equiv 1$ $(\mathrm{mod}\ 2)$. \par
Therefore, by \eqref{31} and \eqref{32}, we obtain the following theorem.
\begin{theorem}
    For $m,p\in\mathbb{N}$ with $m\equiv 1$ $(\mathrm{mod}\ 2)$, we have
    \begin{equation}
    S_{p}^{(k)}(1,m)=\sum_{\nu=0}^{p}\binom{p}{\nu}E_{\nu}^{(k)}\sum_{i=0}^{p-\nu}\binom{p-\nu+1}{i}E_{i}m^{p-i}.\label{33}
    \end{equation}
\end{theorem}
Now, we assume that $p \ge 3$ is an odd integer, so that $E_{p-1}=0$.
Interchanging the order of summation in \eqref{33}, we have
\begin{align}
S_{p}^{(k)}(1,m)\ &=\ \sum_{i=0}^{p}\sum_{\nu=0}^{p-i}\binom{p}{\nu}\binom{p-\nu+1}{i}E_{\nu}^{(k)}E_{i}m^{p-i} \label{34} \\
&=\ \sum_{i=1}^{p-2}\sum_{\nu=0}^{p-i}\binom{p}{\nu}\binom{p-\nu+1}{i}E_{\nu}^{(k)}E_{i}m^{p-i}+\binom{p+1}{p}E_{p}\nonumber \\
&\quad+\sum_{\nu=0}^{p}\binom{p}{\nu}E_{\nu}^{(k)}m^{p}+\sum_{\nu=0}^{1}\binom{p}{\nu}E_{\nu}^{(k)}\binom{p-\nu+1}{p-1}E_{p-1}m\nonumber \\
&=\ \sum_{i=1}^{p-2}\sum_{\nu=0}^{p-i}\binom{p}{\nu}\binom{p-\nu+1}{i}E_{\nu}^{(k)}E_{i}m^{p-i}+(p+1)E_{p}+\sum_{\nu=0}^{p}\binom{p}{\nu}E_{\nu}^{(k)}m^{p}\nonumber \\
&=\ \sum_{i=1}^{p-2}\sum_{\nu=0}^{p-i}\binom{p}{\nu}\binom{p-\nu+1}{i}E_{\nu}^{(k)}E_{i}m^{p-i}+(p+1)E_{p}+m^{p}E_{p}^{(k)}(1). \nonumber
\end{align}
Therefore, we obtain the following theorem.
\begin{theorem}
    For $m\in\mathbb{N}$ with $m\equiv 1$ $(\mathrm{mod}\ 2)$, and $p\equiv 1$ $(\mathrm{mod}\ 2)$ with $p>1$, we have
    \begin{displaymath}
        S_{p}^{(k)}(1,m)=\sum_{i=1}^{p-2}\sum_{\nu=0}^{p-i}\binom{p}{\nu}\binom{p-\nu+1}{i}E_{\nu}^{(k)}E_{i}m^{p-i}+(p+1)E_{p}+m^{p}E_{p}^{(k)}(1).
    \end{displaymath}
In other words, we have
\begin{align}
&m^{p}T_{p}^{(k)}(1,m) \label{35} \\
&\ = \sum_{i=1}^{p-2}\sum_{\nu=0}^{p-i}\binom{p}{\nu}\binom{p-\nu+1}{i}E_{\nu}^{(k)}E_{i}m^{p-i}+(p+1)E_{p}+m^{p}E_{p}^{(k)}(1)+2\sum_{\nu=0}^{p}\binom{p}{\nu}E_{\nu}^{(k)}E_{p+1-\nu}m^{\nu-1}.\nonumber
\end{align}
\end{theorem}
We observe that
\begin{align}
\sum_{\nu=0}^{p}\binom{p-\nu+1}{s}\binom{p}{\nu}E_{\nu}^{(k)}\ &=\ \sum_{\nu=0}^{p-s+1}\binom{p-\nu+1}{s}\binom{p}{\nu}E_{\nu}^{(k)} \label{36} \\
&=\ \sum_{\nu=0}^{p-s}\binom{p-\nu+1}{s}\binom{p}{\nu}E_{\nu}^{(k)}+\binom{p}{s-1}E_{p-s+1}^{(k)}.\nonumber \end{align}
From \eqref{36} and Lemma 8, we have
\begin{align}
&\sum_{\nu=0}^{p-s}\binom{p-\nu+1}{s}\binom{p}{\nu}E_{\nu}^{(k)}    \label{37} \\
&=\sum_{\nu=0}^{p}\binom{p-\nu+1}{s}\binom{p}{\nu}E_{\nu}^{(k)}-\binom{p}{s-1}E_{p-s+1}^{(k)}\nonumber \\
&=\binom{p}{s}E_{p-s}^{(k)}(1)+\binom{p}{s-1}E_{p-s+1}^{(k)}(1)-\binom{p}{s-1}E_{p-s+1}^{(k)}.\nonumber
\end{align}
By \eqref{37}, we get
\begin{align}
&\sum_{i=1}^{p-2}\sum_{\nu=0}^{p-i}\binom{p}{\nu}\binom{p-\nu+1}{i}E_{\nu}^{(k)}E_{i}m^{p-i}\label{38}\\
&=\sum_{i=1}^{p-2}\binom{p}{i}E_{p-i}^{(k)}(1)E_{i}m^{p-i}+\sum_{i=1}^{p-2}\binom{p}{i-1}\big(E_{p-i+1}^{(k)}(1)-E_{p-i+1}\big)E_{i}m^{p-i}.\nonumber
\end{align}
From \eqref{35} and \eqref{38}, we note that
\begin{align}
&m^{p}T_{p}^{(k)}(1,m) \label{39} \\
&=\ \sum_{i=1}^{p-2}\binom{p}{i}E_{p-i}^{(k)}(1)E_{i}m^{p-i}+\sum_{i=1}^{p-2}\binom{p}{i-1}\big(E_{p-i+1}^{(k)}(1)-E_{p-i+1}^{(k)}\big)E_{i}m^{p-i}\nonumber \\
&+(p+1)E_{p}+m^{p}E_{p}^{(k)}(1)+2\sum_{\nu=0}^{p}\binom{p}{\nu}E_{\nu}^{(k)}E_{p+1-\nu}m^{\nu-1}.\nonumber
\end{align}
It is easy to show that
\begin{equation}
E_{1}^{(k)}(1)-E_{1}^{(k)}=1.\label{40}
\end{equation}
By \eqref{39} and \eqref{40}, we get
\begin{align}
&m^{p}T_{p}^{(k)}(1,m)\label{41} \\
&=\sum_{i=0}^{p}\binom{p}{i}E_{p-i}^{(k)}(1)E_{i}m^{p-i}+\sum_{i=1}^{p}\binom{p}{i-1}\big(E_{p-i+1}^{(k)}(1)-E_{p-i+1}^{(k)}\big)m^{p-i}E_{i}\nonumber \\
&\quad+2\sum_{\nu=0}^{p}\binom{p}{\nu}E_{\nu}^{(k)}E_{p+1-\nu}m^{\nu-1}.\nonumber
\end{align}
Therefore, by \eqref{41}, we obtain the following theorem.
\begin{theorem}
For $m\in\mathbb{N}$ with $m\equiv 1$ $(\mathrm{mod}\ 2)$,  $p\equiv 1$ $(\mathrm{mod}\ 2)$ with $p>1$, we have
\begin{align}
&m^{p}T_{p}^{(k)}(1,m)\nonumber \\
&=\sum_{i=0}^{p}\binom{p}{i}E_{p-i}^{(k)}(1)E_{i}m^{p-i}+\sum_{i=1}^{p}\big(E_{p-i+1}^{(k)}(1)-E_{p-i+1}^{(k)}\big)m^{p-i}E_{i} \binom{p}{i-1}\nonumber \\
&\quad+2\sum_{\nu=0}^{p}\binom{p}{\nu}E_{\nu}^{(k)}E_{p+1-\nu}m^{\nu-1}.\nonumber
\end{align}
\end{theorem}
Now, we employ the symbolic notations as $E_{n}(x)=(E+x)^{n}$, $E_{n}^{(k)}(x)=\big(E^{(k)}+x\big)^{n}$, $(n\ge 0)$. \par
Then we first observe that
\begin{align}
&m^{p}\sum_{\mu=0}^{m-1}(-1)^{\mu}\sum_{s=0}^{p}\binom{p}{s}h^{s}E_{s}^{(k)}\bigg(\frac{\mu}{m}\bigg)E_{p-s}\bigg(h-\bigg[\frac{h\mu}{m}\bigg]\bigg) \label{42}\\
&\ = \ m^{p}\sum_{\mu=0}^{m-1}(-1)^{\mu}\bigg(h\bigg(E^{(k)}+\frac{\mu}{m}\bigg)+\bigg(E+h-\bigg[\frac{h\mu}{m}\bigg]\bigg)\bigg)^{p}\nonumber\\
&\ =\ m^{p}\sum_{\mu=0}^{m-1}(-1)^{\mu}\bigg(hE^{(k)}+E+h+\frac{1}{2}-\frac{1}{2}+h\mu m^{-1}-\bigg[\frac{h\mu}{m}\bigg]\bigg)^{p}\nonumber \\
&\ =\ m^{p}\sum_{\mu=0}^{m-1}(-1)^{\mu}\bigg(hE^{(k)}+E+h+\frac{1}{2}+\overline{E}_{1}\bigg(\frac{h\mu}{m}\bigg)\bigg)^{p}.\nonumber
\end{align}
Assume that $h,m$ are relatively prime positive integers. Then, as the index $\mu$ ranges over the values $\mu=0,1,2,\dots,m-1$, the product $h\mu$ does over a complete residue system modulo $m$ and, due to the periodicity of $\overline{E}_{1}(x)$, the term $\overline{E}_{1}\big(\frac{h\mu}{m}\big)$ may be replaced by $\overline{E}_{1}\big(\frac{\mu}{m}\big)$, without alternating the sum over $\mu$.\par
For $m\in\mathbb{N}$ with $m\equiv 1$ $(\mathrm{mod}\ 2)$, and by \eqref{42} and \eqref{18}, we get
\begin{align}
&m^{p}\sum_{\mu=0}^{m-1}(-1)^{\mu}\sum_{s=0}^{p}\binom{p}{s}h^{s}E_{s}^{(k)}\bigg(\frac{\mu}{m}\bigg)E_{p-s}\bigg(h-\bigg[\frac{h\mu}{m}\bigg]\bigg)\label{43}\\
&=\ m^{p}\sum_{\mu=0}^{m-1}(-1)^{\mu}\bigg(hE^{(k)}+h   +E+\frac{1}{2}+\overline{E}_{1}\bigg(\frac{\mu}{m}\bigg)\bigg)^{p}\nonumber \\
&=\ m^{p}\sum_{\mu=0}^{m-1}(-1)^{\mu}\bigg(h(E^{(k)}+1)+E+\frac{\mu}{m}\bigg)^{p}\nonumber \\
&=\ m^{p}\sum_{\mu=0}^{m-1}(-1)^{\mu}\sum_{s=0}^{p}\binom{p}{s}\bigg(E+\frac{\mu}{m}\bigg)^{s}h^{p-s}\big(E^{(k)}+1\big)^{p-s}\nonumber \\
&=\ m^{p}\sum_{\mu=0}^{m-1}(-1)^{\mu}\sum_{s=0}^{p}\binom{p}{s}E_{s}\bigg(\frac{\mu}{m}\bigg)h^{p-s}E_{p-s}^{(k)}(1) \nonumber \\
&=\ \sum_{s=0}^{p}\binom{p}{s}m^{p-s}\bigg(m^{s}\sum_{\mu=0}^{m-1}(-1)^{\mu}E_{s}\bigg(\frac{\mu}{m}\bigg)\bigg)h^{p-s}E_{p-s}^{(k)}(1) \nonumber \\
&=\ \sum_{s=0}^{p}\binom{p}{s}(mh)^{p-s}E_{s}E_{p-s}^{(k)}(1). \nonumber
\end{align}
Therefore, by \eqref{42} and \eqref{43}, we obtain the following theorem.
\begin{theorem}
    For $h,m,p\in\mathbb{N}$ with $(h,m)=1$, $m\equiv 1$ $(\mathrm{mod}\ 2)$, we have
    \begin{align*}
        &m^{p}\sum_{\mu=0}^{m-1}(-1)^{\mu}\sum_{s=0}^{p}\binom{p}{s}h^{s}E_{s}^{(k)}\bigg(\frac{\mu}{m}\bigg)E_{p-s}\bigg(h-\bigg[\frac{h\mu}{m}\bigg]\bigg)\\
        &=\ \sum_{s=0}^{p}\binom{p}{s}(mh)^{p-s}E_{s}E_{p-s}^{(k)}(1).
    \end{align*}
\end{theorem}
From Corollary 7, we note that
\begin{equation}
E_{n}^{(k)}(x)=\sum_{l=0}^{n}\sum_{j=1}^{n+1-l}\sum_{s=0}^{m-1}\binom{n}{l}m^{l}(-1)^{s}E_{l}\bigg(\frac{s+x}{m}\bigg)\frac{S_{1}(n+1-l,j)}{j^{k-1}(n+1-l)},\label{44}
\end{equation}
where $m\in\mathbb{N}$ with $m \equiv 1$ $(\mathrm{mod}\ 2)$, and $n\ge 0$. \par
For $m,h\in\mathbb{N}$ with $m\equiv 1\ (\mathrm{mod}\ 2)$, $h\equiv 1$ $(\mathrm{mod}\ 2)$, by \eqref{44}, we get
\begin{align}
&m^{p}T_{p}^{(k)}(h,m)+h^{p}T_{p}^{(k)}(m,h)\label{45}\\
&=\ 2m^{p}\sum_{\mu=0}^{m-1}\frac{\mu}{m}(-1)^{\mu}\overline{E}_{p}^{(k)}\bigg(\frac{h\mu}{m}\bigg)+2h^{p}\sum_{\nu=0}^{h-1}\frac{\nu}{h}(-1)^{\nu}\overline{E}_{p}^{(k)}\bigg(\frac{m\nu}{h}\bigg)\nonumber \\
&=\ 2m^{p}\sum_{\mu=0}^{m-1}(-1)^{\mu}\frac{\mu}{m}\sum_{l=0}^{p}\binom{p}{l}h^{l}\sum_{\nu=0}^{h-1}\sum_{j=1}^{p+1-l}(-1)^{\nu}\overline{E}_{l}\bigg(\frac{\nu}{h}+\frac{\mu}{m}\bigg)\frac{S_{1}(p+1-l,j)}{j^{k-1}(p+1-l)}\nonumber \\
&+\ 2h^{p}\sum_{\nu=0}^{h-1}(-1)^{\nu}\frac{\nu}{h}\sum_{l=0}^{p}\binom{p}{l}m^{l}\sum_{\mu=0}^{m-1}\sum_{j=1}^{p+1-l}(-1)^{\mu}\overline{E}_{l}\bigg(\frac{\mu}{m}+\frac{\nu}{h}\bigg)\frac{S_{1}(p+1-l,j)}{j^{k-1}(p+1-l)}\nonumber \\
&=\ 2\sum_{\mu=0}^{m-1}(-1)^{\mu}\frac{\mu}{m}\sum_{l=0}^{p}m^{p-l}(mh)^{l}\binom{p}{l}\sum_{\nu=0}^{h-1}\sum_{j=1}^{p+1-l}(-1)^{\nu}\overline{E}_{l}\bigg(\frac{\mu}{m}+\frac{\nu}{h}\bigg)\frac{S_{1}(p-l+1,j)}{(p-l+1)j^{k-1}}\nonumber\\
&+2\sum_{\nu=0}^{h-1}(-1)^{\nu}\frac{\nu}{h}\sum_{l=0}^{p}h^{p-l}(mh)^{l}\binom{p}{l}\sum_{\mu=0}^{m-1}\sum_{j=1}^{p-l+1}(-1)^{\mu}\overline{E}_{l}\bigg(\frac{\nu}{h}+\frac{\mu}{m}\bigg)\frac{S_{1}(p-l+1,j)}{(p-l+1)j^{k-1}}\nonumber \\
&=\ 2\sum_{\mu=0}^{m-1}\sum_{l=0}^{p}\sum_{\nu=0}^{h-1}\sum_{j=1}^{p+1-l}(-1)^{\mu+\nu}(\mu h)(mh)^{-1}m^{p-l}(mh)^{l}\binom{p}{l}\overline{E}_{l}\bigg(\frac{\mu}{m}+\frac{\nu}{h}\bigg)\frac{S_{1}(p-l+1,j)}{(p-l+1)j^{k-1}}\nonumber \\
&+\ 2\sum_{\mu=0}^{m-1}\sum_{l=0}^{p}\sum_{\nu=0}^{h-1}\sum_{j=1}^{p+1-l}(-1)^{\mu+\nu}(\nu m)(mh)^{-1}h^{p-l}(mh)^{l}\binom{p}{l}\overline{E}_{l}\bigg(\frac{\nu}{h}+\frac{\mu}{m}\bigg)\frac{S_{1}(p-l+1,j)}{(p-l+1)j^{k-1}}\nonumber
\end{align}
\begin{displaymath}
=2\sum_{\mu=0}^{m-1}\sum_{l=0}^{p}\sum_{\nu=0}^{h-1}\sum_{j=1}^{p+1-l}(-1)^{\mu+\nu}\frac{(mh)^{l-1}\binom{p}{l}S_{1}(p-l+1,j)}{(p-l+1)j^{k-1}}\big((\mu h)m^{p-l}+(\nu m)h^{p-l}\big)\overline{E}_{l}\bigg(\frac{\nu}{h}+\frac{\mu}{m}\bigg).
\end{displaymath}
Therefore, by \eqref{45}, we obtain the following reciprocity relation.
\begin{theorem}
For $m,h,p\in\mathbb{N}$ with $m\equiv 1$ $(\mathrm{mod}\ 2)$, $h\equiv 1$ $(\mathrm{mod}\ 2)$, and $k\in\mathbb{Z}$,  we have
\begin{align*}
&m^{p}T_{p}^{(k)}(h,m)+h^{p}T_{p}^{(k)}(m,h)\\
& =2\sum_{\mu=0}^{m-1}\sum_{l=0}^{p}\sum_{\nu=0}^{h-1}\sum_{j=1}^{p+1-l}(-1)^{\mu+\nu}\frac{(mh)^{l-1}\binom{p}{l}S_{1}(p-l+1,j)}{(p-l+1)j^{k-1}}\big((\mu h)m^{p-l}+(\nu m)h^{p-l}\big)\overline{E}_{l}\bigg(\frac{\nu}{h}+\frac{\mu}{m}\bigg).
\end{align*}
\end{theorem}
In case of $k=1$, we obtain the following reciprocity relation for the Dedekind type $DC$ sums.

\begin{corollary}
For $m,h,p\in\mathbb{N}$ with $m\equiv 1$ $(\mathrm{mod}\ 2)$, $h\equiv 1$ $(\mathrm{mod}\ 2)$, we have
\begin{align*}
&m^{p}T_{p}(h,m)+h^{p}T_{p}(m,h) \\
&=2(mh)^{p-1}\sum_{\mu=0}^{m-1}\sum_{\nu=0}^{h-1}(-1)^{\mu+\nu}(\mu h+\nu m)\overline{E}_{p}\bigg(\frac{\nu}{h}+\frac{\mu}{m}\bigg).
\end{align*}
\end{corollary}

\section{conclusions}

The generalized Dedekind sums considered by Apostol are given by
\begin{equation*}
S_p(h,m)= \sum_{\mu=1}^{m-1}\frac{\mu}{m}\overline{B}_{p}\big(\frac{h\mu}{m}\big),
\end{equation*}
and satisfy a reciprocity relation (see [1,2]), where $\overline{B}_p(x)$ is the $p$th Bernoulli function.\par
Recently, the type 2 poly-Bernoulli polynomials of index $k$ are defined in terms
of the polyexponential function of index $k$ as
\begin{equation*}
\frac{\mathrm{Ei}_{k}(\log(1+t))}{e^{t}-1}e^{xt}=\sum_{n=0}^{\infty}B_{n}^{(k)}(x)\frac{t^{n}}{n!},\quad(k\in\mathbb{Z}).
\end{equation*} \par
As a further extension of the generalized Dedekind sums, the poly-Dedekind sums defined by
\begin{equation*}
S_p^{(k)}(h,m)= \sum_{\mu=1}^{m-1}\frac{\mu}{m}\overline{B}_{p}^{(k)}\bigg(\frac{h\mu}{m}\bigg),
\end{equation*}
were considered and shown to satisfy a reciprocity relation in
[16], where $\overline{B}_p^{(k)}(x)=B_{p}^{(k)}(x-[x])$ are the
type 2 poly-Bernoulli functions of index $k$, and
$S_p^{(1)}(h,m)=S_p(h,m)$. \par The Dedekind type $DC$ sums
defined by
\begin{equation*}
T_{p}(h,m)=2\sum_{\mu=1}^{m-1}(-1)^{\mu}\frac{\mu}{m}\overline{E}_{p}\bigg(\frac{h\mu}{m}\bigg)
\end{equation*}
were introduced and shown to satisfy a reciprocity relation in
[14], where $\overline{E}_{p}(x)$ is the $p$-th Euler function.
Simsek found trigonometric representations of the Dedekind type DC
sums and their relations to  Clausen functions, polylogarithm
function, Hurwitz zeta function, generalized Lambert series
(G-series), and Hardy-–Berndt sums. \par In this paper, as a
furthrer generalization of the Dedekind type DC sums, the
poly-Dedekind type DC sums given by
\begin{equation*}
T_{p}^{(k)}(h,m)=2\sum_{\mu=1}^{m-1}(-1)^{\mu}\frac{\mu}{m}\overline{E}_{p}^{(k)}\bigg(\frac{h\mu}{m}\bigg)
\end{equation*}
were considered and shown, among other things, to satisfy a reciprocity relation in Theorem 14. \par
Finally, we observe that the Dedekind sums and their generalizations are defined in terms of Bernoulli functions and their generalizations, while the Dedekind type DC sums and their generalizations are defined in terms of Euler functions and their generalizations.

\vspace{0.2in}

\noindent{\bf{Acknowledgments}} \\
       We thank Jangjeon Research Institute for Mathematical Sciences for its support during the preparation of this paper.

\vspace{0.2in}

\noindent{\bf{Conflict of interest}} \\
All authors declare no conflicts of interest in this paper.

\end{document}